\documentclass[12pt,a4paper]{amsart}
\usepackage{graphicx}
\usepackage{epstopdf}
\usepackage{amssymb}
\usepackage{amsthm}
\usepackage{amsmath}
\usepackage{pstcol, pst-node}
\usepackage{mathrsfs}
\usepackage{amscd}

\newtheorem{thm}{Theorem}[section]

\newtheorem{assum}[thm]{Assumption}
\newtheorem{cor}[thm]{Corollary}
\newtheorem{prop}[thm]{Proposition}
\newtheorem{lem}[thm]{Lemma}
\newtheorem{rem}[thm]{Remark}
\newtheorem{claim}[thm]{Claim}

\theoremstyle{definition}
\newtheorem{defn}[thm]{Definition}

\newtheorem{prop-def}[thm]{Proposition-Definition}

\theoremstyle{remark}

\newcommand{\be}{\begin{equation}}
\newcommand{\bc}{\begin{cor}}
\newcommand{\bt}{\begin{thm}}
\newcommand{\bl}{\begin{lem}}
\newcommand{\bpr}{\begin{prop}}
\newcommand{\br}{\begin{rem}}
\newcommand{\bd}{\begin{defn}}
\newcommand{\ee}{\end{equation}}
\newcommand{\et}{\end{thm}}
\newcommand{\el}{\end{lem}}
\newcommand{\epr}{\end{prop}}
\newcommand{\er}{\end{rem}}
\newcommand{\ed}{\end{defn}}
\newcommand{\ec}{\end{cor}}
\newcommand{\R}{\Bbb{R}}

\newcommand{\nnn}{\noindent}

\newcommand{\del}{\partial}

\newcommand{\Z}{\Bbb{Z}}

\begin{document}
\title{Integration of vector fields on cell complexes and Morse theory} 
\author{Takeo Nishinou}
\address{Department of Mathematics, Rikkyo University,
  Toshima, Tokyo, Japan } 
\email{nishinou@rikkyo.ac.jp}
\begin{abstract}

In this paper, we investigate vector fields on polyhedral complexes and
 their associated trajectories.
We study vector fields which are analogue of the gradient vector field of a function in the smooth case.
Our goal is to define a nice theory of trajectories of such vector fields, 
 so that the set of them captures the topology
 of the polyhedral complex, as in classical Morse theory.
Since we do not assume the polyhedral complex to be a manifold, the definition of vector fields
 on it is very different from the smooth case.
Nevertheless, we will show that there exist nice classes of functions and metrics which give
 gradient vector fields with desired properties.
Our construction relies on Forman's discrete Morse theory.
In particular, the class of functions we use is an improvement of 
 Forman's discrete Morse functions.
A notable feature of our theory is that 
 our gradient vector fields are defined purely from functions
 and metrics as in the smooth case, contrary to the case of discrete Morse theory  
 where we need the data of dimension of cells.
This allows us to implement several useful constructions which were 
 not available in the discrete case.

\end{abstract}
\maketitle
\section{Introduction}
In this paper, we consider vector fields on cell complexes.
Our goal is to give a reasonable theory of integration of such vector fields.
In other words, we would like to give a definition of trajectories of these vector fields
 with nice properties.
As a measure of such `nice properties', we refer to Morse theory \cite{M}.
Namely, if we can calculate the homology of cell complexes in a simple way
 from the data of trajectories
 of a vector field associated with a suitable function, it may be reasonable to claim that the theory of integrating vector fields is nice.
 
Since we deal with cell complexes which are not necessarily manifolds, a definition of vector fields
 should be different from the usual one.
In particular, to capture the topological informations of a cell complex, it will be necessary to allow
 multivaluedness at points where there is no regular neighborhood.
On the other hand, if we allow any object with such multivaluedness, 
 we will have to deal with extremely complicated objects which exhibit pathological behavior.
Therefore, it is desirable to assume that vector fields are as simple as possible on each cell.

From this point of view,  
 we consider gradient vector fields of piecewise linear functions.
However, it is easy to see that if one chooses the functions arbitrarily, 
 it would be impossible to deduce any meaningful result, even if one chooses them generically.
Different choices of functions will give different results, and there would be no reasonable control.
Therefore, for the construction of a theory of integration of gradient vector fields, 
 it is essential to introduce a nice class of functions.

There are discrete, or piecewise linear, analogues of Morse theory, see for example \cite{BB, F1, Z}.
We adopt Forman's \emph{discrete Morse functions} \cite{F1}, see Definition \ref{def:dmorsefunc},
 to solve the problem.
Forman used such functions to develop \emph{discrete Morse theory}, a nice combinatorial analogue
 of usual Morse theory on smooth manifolds.
However, if we allow arbitrary discrete Morse functions, we will still have difficulties.
We need to choose a nice class of discrete Morse functions, which we call \emph{tame},
 to have a good control of gradient trajectories, see Definition \ref{def:morsefunctionprop}.

On the other hand, to consider the gradient vector field of a function, we also need a metric.
It turns out that the choice of a suitable class of metrics is also essential.
We choose \emph{sharp 
 piecewise affine metrics} on the base space, see Definition \ref{def:sharpmetric}.
This is also necessary to guarantee nice behavior of trajectories of vector fields.
On the other hand, on any finite cell complex, plenty of tame discrete Morse functions and sharp
 piecewise affine metrics
 exist, so these choices do not give a restriction to the theory.

Using a tame discrete Morse function and a sharp piecewise affine metric, we can define the gradient vector field
 and its trajectories, see Definition \ref{def:gradfield} and \ref{def:gradtraj}.
They are the analogues of the corresponding objects in smooth case.
However, the nature of gradient vector fields defined here is quite different from 
 those in the smooth case.
For example, even when the base space is a manifold, 
 several flows can go into (or emanate from)
 a point, which is not necessarily a critical point, see, for example, Figure \ref{fig:1}.
Nevertheless, the total behavior of the flows is good enough and we can 
 naturally define a complex from these flows which computes the homology
 of the base space, see Section \ref{sec:8}.
We will prove this by comparing it to discrete Morse theory, see Theorem \ref{thm:main}.

As we mentioned at the beginning, our study is an attempt to construct a theory of vector fields on
 cell complexes in a natural and tractable way.
This is achieved by geometrically realizing discrete Morse theory.
One of the nice features of our construction is that the gradient vector field of a function
 is defined purely from the information of the function and the metric, as in the smooth case.
This is different from the case of
 discrete Morse theory, where the definition of gradient paths requires 
 information other than the function and the metric, namely the dimension of the cell.

This allows us to give another useful construction which was not available in discrete Morse theory.
Namely, when the base space is a closed manifold, we can consider the dual theory.
In discrete Morse theory, the dual theory is defined on the dual cell complex, 
 and one cannot directly compare the gradient vector fields of the original theory with those of the dual theory.
On the other hand, in our case, the gradient vector field of the dual theory is obtained
 simply by reversing the directions, as in the smooth case.
This allows us to construct stable and unstable complexes.
In fact, we can introduce stable and unstable complexes even when the base space is not a manifold.
In such a case, the reversed vector field is not necessarily related to a discrete Morse theory, since
 the dual complex of a polyhedral complex may not make sense when the base space is not a manifold.
Nevertheless, the vector field still has nice properties, 
 see Section \ref{sec:stablecomplex}.


\section{Review of discrete Morse theory}\label{sec:discreteMorse}

Our purpose is to study the behavior of gradient vector fields associated with piecewise linear functions.
To make sense of such vector fields, the base space needs to have a
 piecewise linear structure.
Simplicial complexes  or triangulated manifolds naturally have such structures.
We will work with affine complexes which we now introduce.

\begin{defn}\label{def:polytope}
A \emph{convex polytope} $P$ is a convex hull of finite points in 
 the vector space $\R^n$ for some $n$.
The boundary $\partial P$ is naturally a union $\cup P_i$ of convex polytopes 
 such that for any $i\neq j$, $P_i\cap P_j = P_k$ for some $k$, if $P_i\cap P_j$ is not empty.
Each $P_i$ is called a \emph{face} of $P$.
A face of codimension one in $P$ is called a \emph{facet}.
Note that any polytope has a natural affine linear structure.
That is, the notion of affine linear functions on a polytope makes sense.
\end{defn}

\begin{defn}
An \emph{affine complex} $X$ is a finite CW complex such that
 each cell has a structure of a convex polytope in the sense of 
 Definition \ref{def:polytope} compatible with the gluing.
Here, a cell always means a closed cell in this paper.
Compatibility means that if $X_j$ is a face of some cell $X_i$ of $X$, 
 the affine linear structure on $X_j$ induced by its structure as a convex polytope
 is isomorphic to the affine linear structure induced by that of $X_i$.
Note that no self-intersection of a cell is allowed.
\end{defn}
In particular, an affine complex need not be a topological manifold,
 and it can have boundary components.
\begin{defn}
The \emph{dimension} of an affine complex $X$ is the maximum of the dimension of 
 the cells in $X$.
\end{defn}
Let $X$ be an affine complex.
If $\alpha$ is an $i$-dimensional cell of $X$, we write it as $\alpha^{(i)}$ when we want to 
 emphasize its dimension.
We often call it an $i$-cell.
If $\alpha^{(i)}$ is a face of $\beta^{(j)}$, we write
 $\alpha^{(i)} < \beta^{(j)}$.
In this case, we say that $\alpha^{(i)}$ is adjacent to $\beta^{(j)}$, and also that
 $\beta^{(j)}$ is adjacent to $\alpha^{(i)}$.
 
As we mentioned in the introduction, we study vector fields modeled on
 the gradient vector fields of piecewise linear functions.
To introduce a nice class of functions which are compatible with the piecewise linearity, we use {discrete
 Morse functions} in the sense of Forman \cite{F1}, 
 instead of smooth functions in usual Morse theory.
Let us recall the definition.


%
%
\begin{defn}\label{def:dmorsefunc}(\cite[Definition 2.1]{F1})
A real valued function
\begin{equation*}
F\colon \{ \rm{cells \; of \; X} \} \to \R
\end{equation*}
is a \emph{discrete Morse function} if for every cell $\alpha^{(i)}$ of $X$,
 the inequalities\\
\begin{equation}\label{eq:morse1}
\# \{\beta^{(i+1)} > \alpha^{(i)}\; |\; F(\beta^{(i+1)}) \leq F(\alpha^{(i)})\} \leq 1
\end{equation}
and 
\begin{equation}\label{eq:morse2}
\# \{\gamma^{(i-1)} < \alpha^{(i)} \;|\; F(\gamma^{(i-1)}) \geq F(\alpha^{(i)})\} \leq 1
\end{equation}
 hold.
\end{defn}
Roughly speaking, 
 this definition says a discrete Morse function tends to have larger values
 on larger dimensional cells.
For example, the function $F_{triv}$ defined by
\[
F_{triv}(\alpha^{(i)}) = i
\]
 for any $i$-dimensional cell $\alpha^{(i)}$ is a discrete Morse function,
 where the left hand side of 
 (\ref{eq:morse1}) and (\ref{eq:morse2}) are both zero.
The condition of the definition claims that there is at most one
 exception for each cell.
The cells without such an exception is called critical as in the following definition.
\begin{defn}(\cite[Definition 2.2]{F1})
Consider an affine complex and a discrete Morse function $F$ on
 it.
A cell $\alpha^{(i)}$ is \emph{critical} if
\begin{equation}
 \#\{\beta^{(i+1)} > \alpha^{(i)}\; |\; F(\beta^{(i+1)}) \leq F(\alpha^{(i)})\} = 0
\end{equation}
and
\begin{equation}
 \# \{\gamma^{(i-1)} < \alpha^{(i)}\; |\; F(\gamma^{(i-1)}) \geq F(\alpha^{(i)})\} = 0
\end{equation}
 hold.
\end{defn}
For example, when we take the discrete Morse function
 $F_{triv}$ above, all the cells are critical.

We also need the definition of gradient flows in 
 discrete Morse theory, which is used to define the
 discrete Morse complex.
We need it because later we will compare
 the discrete Morse complex with a complex defined from
 piecewise linear gradient flows, see Section \ref{sec:9}.

Recall that there is a chain complex naturally associated with a CW complex.
Let $C_i(X, \Z)$ denote the free abelian group generated 
 by $i$-dimensional cells of $X$,
 where on each cell an orientation is chosen.
In particular, $-\sigma \in C_i(X, \Z)$
 represents the same cell $\sigma$ with the opposite orientation.
Let $\del$ denote the boundary operator
\begin{equation*}
\del : C_{i+1}(X, \Z) \to C_{i}(X, \Z).
\end{equation*}
For example, in the simplicial case,
 it is given by
\begin{equation*}
\del \langle a_0, \dots, a_{i+1}\rangle = \sum_{l=0}^{i+1} (-1)^l
    \langle a_0, \dots, \check{a}_l, \dots, a_{i+1}\rangle,
\end{equation*}
 here $\check{a}_l$ means it is skipped.
In general, we have
\[
\partial\beta = \sum_{l}sgn(\beta, \alpha_l)\alpha_l,
\]
 where the sum runs over the set of facets of $\beta$, 
 and $sgn(\beta, \alpha_l)\in\{\pm 1\}$ is determined as follows.
Namely, given an oriented polytope $\beta^{(i+1)}$, the orientation is represented by 
 a non-zero element $\tau_{i+1, x}\in \wedge^{i+1} T_x\beta^{(i+1)}$ at any point $x$ in the interior of $\beta^{(i+1)}$.
This determines a constant section $\tau_{i+1}$ of $\wedge^{i+1} T\beta^{(i+1)}$.
On a facet $\alpha^{(i)}$ of $\beta^{(i+1)}$, 
 there is a natural orientation represented by an element $\eta_i\in \wedge^i T_y\alpha^{(i)}$
 at any point $y$ in the interior of $\alpha^{(i)}$, determined by the condition 
 $n_y\wedge \eta_i = \tau_{i+1, y}$, where $n_y$ is an outward normal vector of $\beta^{(i+1)}$ at $y$.
Then,
\[
sgn(\beta, \alpha_l)
 = \begin{cases}
   1\;\; \text{if the orientation on $\alpha_l$ induced from $\beta$ equals to the fixed orientation on $\alpha_l$}, \\
   -1\;\; {\rm otherwise}.
 \end{cases}
\]

Now, introduce an inner product $\langle\, , \,\rangle$ on $C_*(X, \Z)$
 by requiring the set of cells of $X$ to be an orthonormal basis.
\begin{defn}(\cite[Definition 6.1]{F1},)
A \emph{gradient vector field} $V$ on an affine complex equipped with
 a discrete Morse function $F$ and an orientation for
 each cell is defined as follows.
Let $\alpha^{(i)}$ be an oriented $i$-cell.
If there is an $(i+1)$-cell $\beta^{(i+1)}$
 satisfying $\beta^{(i+1)} > \alpha^{(i)}$ and $F(\beta^{(i+1)}) \leq F(\alpha^{(i)})$,
 we set
\begin{equation*}
V(\alpha^{(i)}) = -\langle\del \beta^{(i+1)}, \alpha^{(i)}\rangle \beta^{(i+1)}.
\end{equation*}
If there is no such $\beta^{(i+1)}$, set
\begin{equation*}
V(\alpha^{(i)}) = 0.
\end{equation*}
For each $i$, we extend $V$ linearly to a map
\begin{equation*}
V: C_i(X, \Z) \to C_{i+1}(X, \Z).
\end{equation*}
The map $V$ is a discrete analogue of the gradient vector field of a 
 Morse function on a smooth manifold.
The next is an analogue of the integral of such a vector field.
\end{defn}
\begin{defn}\label{def:gradpath}(\cite[Definition 8.4]{F1})
Let $X$ be an affine complex equipped with a 
 discrete Morse function $F$.
A $gradient \,\, path$ of dimension $i$ is a sequence $\gamma$ of $i$-cells
 of $X$
\begin{equation*}
\gamma = \alpha_0, \alpha_1, \alpha_2, \dots, \alpha_r
\end{equation*}
such that for every $l = 0, \dots, r-1$,\\

(i) if $V(\alpha_l) = 0$, $\alpha_{l+1} = \alpha_l$.\\

(ii) if $V(\alpha_l) \neq 0$, $\alpha_{l+1} < V(\alpha_l)$
 and $\alpha_{l+1} \neq \alpha_l$.\\

We call $\gamma$ a gradient path of $length$ $r$.

\end{defn}

\subsection{Sign of gradient paths and Morse differential}\label{subsec:sign1}

Let us introduce the complex of Forman
\begin{equation*}
(\ast) \;\;\;\;
\mathcal M: 0 \to \mathcal M_n \xrightarrow{\tilde{\del}}
 \mathcal M_{n-1} \xrightarrow{\tilde{\del}} \cdots
 \xrightarrow{\tilde{\del}} \mathcal M_{0} \to 0, 
\end{equation*}
in terms of gradient paths, as in \cite[Section 8]{F1}.
Here, $\mathcal M_i$ is the free abelian group generated by
 critical $i$-dimensional cells.

For this purpose, we need a remark on the orientation
 of the cells contained in a gradient path, as in \cite[page 125]{F1}.
Suppose $\alpha$ and $\tilde \alpha$ are distinct $i$-cells
 and $\beta$ is an $(i+1)$-cell with $\alpha < \beta$
 and $\tilde \alpha < \beta$.
Then, an orientation on $\alpha$ induces an orientation on $\beta$
 by the condition $\langle\del \beta, \alpha\rangle = -1$. 
Given this orientation on $\beta$, we choose the orientation on 
 $\tilde \alpha$ so that $\langle\del \beta, \tilde \alpha\rangle = 1$.
Equivalently, fixing an orientation on $\alpha$ and $\beta$,
 an orientation is induced on $\tilde \alpha$ by the condition
\begin{equation*}
\langle\del \beta, \alpha\rangle\langle\del \beta, \tilde \alpha\rangle = -1.
\end{equation*}
On the other hand, if $\alpha = \tilde \alpha$, we just take the same orientation on $\tilde\alpha$
 as $\alpha$.
Thus, if $\gamma = \alpha_0, \alpha_1, \dots, \alpha_r$
 is a gradient path, an orientation on $\alpha_0$ induces
 an orientation on each $\alpha_i$, and, in particular, on $\alpha_r$.
Recall that we have fixed an orientation for each cell of $X$.
Write $m(\gamma) = 1$ if the fixed orientation on $\alpha_0$
 induces the fixed orientation on $\alpha_r$, and
 $m(\gamma) = -1$ otherwise. 
For $i$-cells $\alpha$ and $\tilde \alpha$, let
 $\Gamma_r(\alpha, \tilde \alpha)$ denote the set of all gradient
 paths from $\alpha$ to $\tilde \alpha$ of length $r$.
If $\beta^{(i+1)}$ and $\alpha^{(i)}$ are critical,
 the differential $\tilde \del$ of the sequence ($\ast$) is given by
\begin{equation}\label{eq:discretediff}
\langle\tilde \del \beta^{(i+1)}, \alpha^{(i)}\rangle =
 \sum_{\tilde \alpha^{(i)} < \beta^{(i+1)}}\langle\del \beta^{(i+1)}, \tilde \alpha^{(i)}\rangle
 \sum_{\gamma \in \Gamma_N(\tilde \alpha^{(i)}, \alpha^{(i)})} m(\gamma)
\end{equation}
for any $N$ large enough.
In the paper \cite{F1}, Forman proved the next theorem
 (for any finite CW complex which is not necessarily an affine complex).
\begin{thm}\label{thm:discmtheory}(Forman \cite{F1}, Theorems 8.2 and 8.10)
$\tilde \del$ is a differential, i.e, $\tilde \del^2 = 0$.
The homology of the complex $(\ast)$
 is precisely the singular homology of $X$. \qed
\end{thm}

\section{Tame discrete Morse functions}
Our purpose is to define and study integration of  vector fields on cell complexes
 using discrete Morse theory as a model.
However, if one uses arbitrary discrete Morse functions, one will soon realize that
 the corresponding piecewise linear theory may be ill-behaved.
To remedy this, we need several ideas, as we explained in the introduction.
The first of these ideas is to choose a nice class of discrete Morse functions, which we call tame.

\begin{defn}\label{def:morsefunctionprop}
We call a discrete Morse function $F$ \emph{generic} if for any pair of cells $\alpha, \beta$
 satisfying $\alpha<\beta$, we have
 $F(\alpha)\neq F(\beta)$.
We call $F$ \emph{tame} if $\alpha^{(i)} < \beta^{(j)}$ and
 $j\geq i+2$, the inequality $F(\beta^{(j)}) > F(\alpha^{(i)})$ holds.
\end{defn}

The tameness of discrete Morse functions does not impose a strong restriction.
In fact, given a generic discrete Morse function, we can modify it to a tame one 
 without changing discrete Morse theoretic properties.
\begin{lem}\label{lem:tamefunction}
Let $F$ be a generic discrete Morse function on an affine complex $X$.
Then, there is a tame generic discrete Morse function $\tilde F$ such that there is a
 natural identification between the 
 associated complexes $\mathcal M$.
That is, the set of critical cells and the differential on them are the same.
\end{lem}
\proof
We prove this by induction on the number of pairs of cells $(\alpha^{(i)}, \beta^{(j)})$
 in $X$
 which violate the condition for the tameness.
That is, we consider those pairs of cells which satisfy
\[
\alpha^{(i)} < \beta^{(j)}, \;\; j \geq i+2,\;\; F(\alpha^{(i)}) > F(\beta^{(j)}). 
\]
If this number is zero, $F$ itself is tame.
Let $k$ be a positive integer.
Assume that for any $F$ such that there are $k'$-pairs of cells $(\alpha, \beta)$, 
 $k'\in \{1, 2, \dots, k\}$,
 which violate the tameness condition, there is a tame $\tilde F$ which gives the same
 complex as the complex associated with $F$.
We will prove the assertion for $F$ for which there are $(k+1)$-pairs of cells $(\alpha, \beta)$
 which violate the tameness condition.

Let $\mathcal P$ be the set of such pairs of cells.
We introduce a partial order to $\mathcal P$
 by the rule that
 $(\alpha^{(i)}, \beta^{(j)}) > ((\alpha')^{(i')}, (\beta')^{(j')})$ if and only if 
 $j> j'$.
Let $(\alpha^{(i)}, \beta^{(j)})$ be a maximal element with respect to this order.
Note that if $\alpha^{(i)} < \beta^{(j)}$ and $j\geq i+2$, 
 there are at least two $(i+1)$-cells adjacent to both 
 $\alpha^{(i)}$ and $\beta^{(j)}$. 
Then, by definition of discrete Morse functions, 
 there is at least one $(i+1)$-cell $\gamma^{(i+1)}$ such that
\[
\alpha^{(i)} < \gamma^{(i+1)} < \beta^{(j)}
\]
 and
\[
F(\gamma^{(i+1)})> F(\alpha^{(i)})
\]
 hold.
Similarly, if we have $j\geq i+3$, 
 there is at least one $(i+2)$-cell $\delta^{(i+2)}$ such that 
\[
\gamma^{(i+1)} < \delta^{(i+2)} < \beta^{(j)}
\]
 and
\[
F(\delta^{(i+2)})> F(\gamma^{(i+1)})
\]
 hold.
In particular, we have $\alpha^{(i)} < \delta^{(i+2)}$ and $F(\delta^{(i+2)}) > F(\alpha^{(i)})$.
Repeating this, we see that there is at least one $(j-1)$-cell $\varepsilon^{(j-1)}$ 
 which satisfies
\[
\alpha^{(i)} < \varepsilon^{(j-1)} < \beta^{(j)}
\]
 and 
\[
F(\varepsilon^{(j-1)}) > F(\alpha^{(i)}) > F(\beta^{(j)}).
\]
In particular, all the $(j-1)$-cells other than $\varepsilon^{(j-1)}$ adjacent to $\beta^{(j)}$
 have the smaller value of $F$ than $\beta^{(j)}$ by definition of discrete Morse functions.

On the other hand, if $\eta^{(k)} > \beta^{(j)}$
 (so that we also have $\eta^{(k)} >\varepsilon^{(j-1)}$), 
 by the maximality of the pair $(\alpha^{(i)}, \beta^{(j)})$
 with respect to the given order on $\mathcal P$, we have
\begin{equation}\label{eq:-1}
F(\eta^{(k)}) > F(\varepsilon^{(j-1)}).
\end{equation}

Then, modify the value $F(\beta^{(j)})$ to $\bar F(\beta^{(j)})$,
 where $\bar F(\beta^{(j)})$ is a generic real number in the open interval 
$
(F(\alpha^{(i)}), F(\varepsilon^{(j-1)})).
$
The values of $\bar F$ at the other cells are the same as those of $F$.
Clearly, $\bar F$ is still a generic discrete Morse function.
By the relation (\ref{eq:-1}),
 this modification does not produce a new pair of cells which violates
 the tameness condition.
Also, it does not change the associated complex.
Then, by applying the induction hypothesis to $\bar F$, we obtain a required function $\tilde F$.\qed

\section{Piecewise linear functions associated with discrete Morse functions
 and metrics on affine complexes}
Now, we construct a piecewise linear function $f$ from 
 a given discrete Morse function $F$.
\begin{defn}\label{def:morsefunc}
Let $X$ be an affine complex and $F$ be a discrete
 Morse function on it.\\

\noindent
1. Take the barycentric subdivision $X_1$ of $X$. By definition, the vertices of it
 correspond to the cells of $X$.\\

\nnn
2. To a vertex of $X_1$, give the value of $F$ at
 the corresponding
 cell.\\

\nnn
3. Extend it affine linearly to each maximal cell of $X_1$.
 It is possible, because $X_1$ is a simplicial complex.  \\

We call a vertex of $X_1$ corresponding to an $i$-dimensional cell
 of $X$ an $i$-dimensional vertex.
Also, we call an edge of $X_1$ connecting an $i$-dimensional vertex
 and a $j$-dimensional vertex  an $(i,j)$-edge.
\end{defn}

One of the features of discrete Morse theory is 
 that although it studies an analogue of gradient flows, 
 it does not require a metric on the base space.
However, 
 since we consider a geometric flow on $X$ (contrary to the
 combinatorial one, see Definition \ref{def:gradpath}),
 we need a metric as in usual Morse theory.

\begin{defn}
Let $Y$ be a simplicial complex.
A \emph{piecewise affine metric} on $Y$ is defined by  
 an affine metric on each maximal cell 
 (i.e, a cell which is not contained in the boundary of
 another cell) which is compatible in the following sense.
Namely, on any lower dimensional cell, the metrics induced
 by those on the adjacent maximal cells are the same.
Here, an affine metric on a polytope in $\Bbb R^n$ is the metric obtained by the restriction
 of a constant valued positive definite symmetric 2-tensor field $g\in \Gamma(\Bbb R^n, Sym^2T^*\Bbb R^n)$.
\end{defn}  

We cannot use every such a metric, since otherwise 
 the gradient flow will exhibit 
 undesirable behavior. 
This problem can be resolved by imposing the following simple
 condition to the metric.
\begin{defn}\label{def:sharpmetric}
We call a piecewise affine metric on a simplicial complex
 \emph{sharp} if on each simplex $\sigma$, the following condition is satisfied.
Namely, take any vertex $a$ of $\sigma$, 
 and a two dimensional plane $L$ containing $a$.
Assume the intersection $L\cap \sigma$ is a triangle, which we write by $\triangle abc$.
Then, it is acute-angled, that is, the interior angle at each vertex is less than $\frac{\pi}{2}$
 with respect to the induced metric.
%
\end{defn}
A sharp piecewise linear metric always exists on simplicial complexes as the following claim shows.
\begin{prop}\label{prop:sharpmetric}
Given a simplicial complex, a piecewise affine metric for which 
 any of its simplex is equilateral of a given edge length is sharp.
\end{prop}
\proof
It suffices to prove that the metric on the standard simplex with vertices
 at $a_0 = (1, 0, \dots, 0), a_1 = (0, 1, 0, \dots, 0), \dots, a_n = (0, \dots, 0, 1)$
 in $\Bbb R^{n+1}$ induced from the Euclidean metric on $\Bbb R^{n+1}$ is sharp.
By symmetry, we can assume that the vertex $a$ in Definition \ref{def:sharpmetric}
 is $a_0$.
Then, using the notation in Definition \ref{def:sharpmetric}, 
 we can assume that the vertices $b$ and $c$ have the coordinates of the form
\[
(0, x_1, \dots, x_k, 0, \dots, 0), \;\; x_1, \dots, x_k\geq 0,
\]
 and
\[
(0, \dots, 0, y_{k+1}, \dots, y_n), \;\; y_{k+1}, \dots, y_n\geq 0,
\]
 respectively.
Here, the coordinates satisfy $\sum_{i=1}^kx_i = 1$ and $\sum_{i=k+1}^ny_i = 1$.
 
By the cosine formula, we have
\[
\cos\angle acb = \frac{\overline{ac}^2 + \overline{bc}^2 - \overline{ab}^2}{2\overline{ac}\overline{bc}}.
\]
Since we have
\[
\overline{ac}^2 + \overline{bc}^2 - \overline{ab}^2
 = 2\sum_{i=k+1}^n y_i^2 > 0,
\]
 we have $\cos\angle acb > 0$ and $\angle acb < \frac{\pi}{2}$.
We can apply the same argument to the other angles.
This finishes the proof.\qed

\section{Piecewise linear gradient vector fields}\label{sec:5}
In Sections \ref{sec:5} to \ref{sec:8}, we will define piecewise linear Morse theory on affine complexes.
First, we will define the piecewise linear gradient flow.
Let $X, X_1, F$ and $f$ be as in the previous section.
We write $n = \dim X$.
We fix a piecewise affine metric on $X_1$.
From now on, we assume the following unless otherwise noted.

\begin{assum}\label{assum:1}
The metric is sharp and 
 the discrete Morse function $F$ is generic and tame in the sense of Definition \ref{def:morsefunctionprop}.
\end{assum}
By Lemma \ref{lem:tamefunction} and Proposition \ref{prop:sharpmetric}, 
 this assumption does not give a restriction to the applicability of the theory.
 
Let $\sigma$ be a cell in $X_1$.
If $\sigma$ is a polytope in $\Bbb R^n$ and $L$ is the minimal affine subspace of $\Bbb R^n$ containing 
 $\sigma$, we define
\[
T\sigma = TL|_{\sigma}.
\]
The restriction of the function $f$ and the metric to $\sigma$ induces
 a gradient vector field on it.
By the piecewise linearity, this is a constant vector field.
Under the above assumption, we have the following.
\begin{lem}\label{lem:generaldirection}
Under Assumption \ref{assum:1}, the direction of the gradient vector on each $\sigma$ is not 
 contained in the hyperplanes
 parallel to the facets of it.\qed
\end{lem}

\begin{defn}
Let $\sigma$ be a cell in $X_1$.
The constant vector field above gives
 the \emph{gradient flow} on $\sigma$, which is the family of maps 
\[
\phi_t\colon \sigma\to \sigma,\;\; t\in \Bbb R_{\geq 0}.
\]
Note that since $\sigma$ has boundary, $\phi_t$ is not a diffeomorphism.
Namely, once a point reaches the boundary by the flow at some $t>0$, it stays there afterwards.
\end{defn}

\begin{defn}\label{def:in-out}
Let $\sigma$ be a cell in $X_1$ and $\tau$ be its facet.
Let $x$ be a point in the interior of $\tau$.
By Lemma \ref{lem:generaldirection}, for sufficiently small $t$, 
 either
\begin{enumerate} 
\item[(i)] $\phi_t(x)$ is in the interior of $\sigma$, or 
\item[(ii)] $\phi_t(x) = x$,
\end{enumerate}
 holds.
In the case (i), we say that the facet $\tau$ has an \emph{out-flow},
 and in the case (ii), we say that the facet $\tau$ has an \emph{in-flow}.
\end{defn}

Let $x$ be a point in $X_1$.
Let $\sigma_0$ be the minimal cell of $X_1$ which contains $x$ in its interior.
We write its dimension by $d_0 = \dim \sigma_0$.
Let $\{\sigma_{j, 1}, \dots, \sigma_{j, a_j}\}$ be the set of $(d_0+j)$-dimensional cells
 which contain $\sigma_0$, here $1\leq j\leq n-d_0$.
For each $\sigma_{j, l}$, there is a gradient vector field determined by the restriction of $f$
 and the metric.
In this paper, we assume that the gradient vector points the direction in which 
 the function decreases.
We write the value of it at $x$ as $\tau_{x;j, l}\in T_{x}\sigma_{j, l}$.
Here, we write $\sigma_0 = \sigma_{0, 1}$.
The inclusions of the cells induce natural inclusions of tangent spaces.
Note that under the above assumption, the images of gradient vectors 
 under these inclusions are all different.

\begin{figure}[h]
\includegraphics{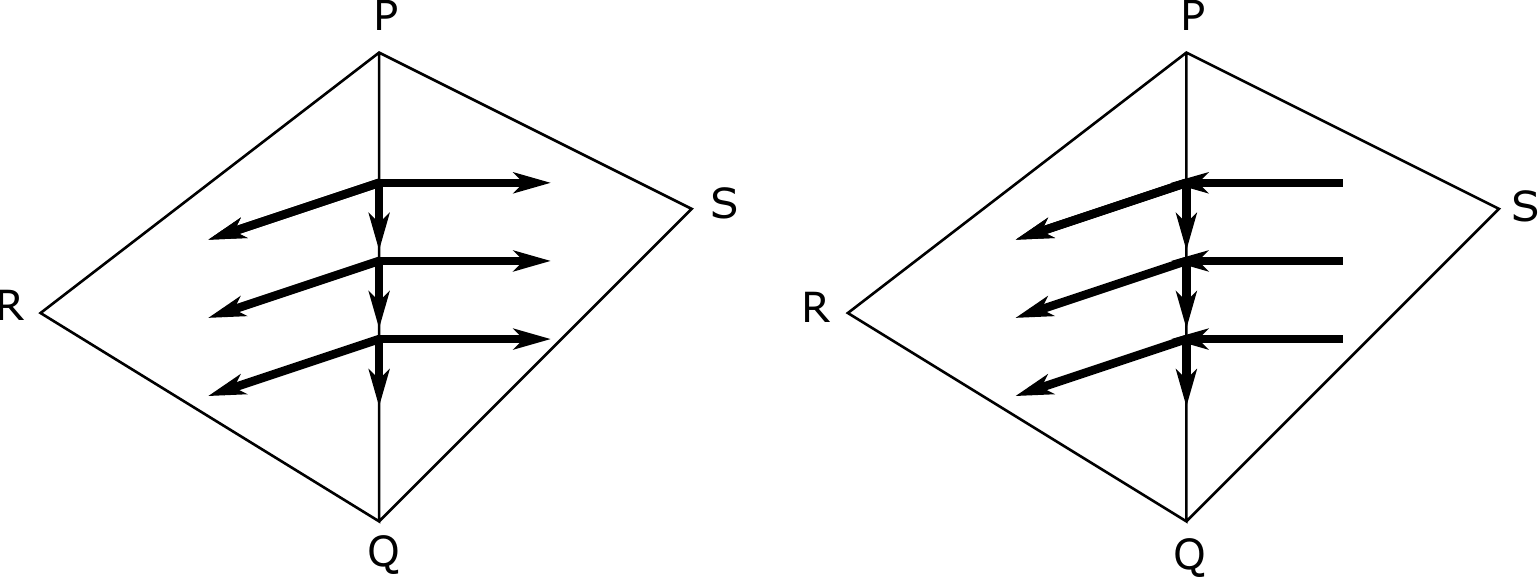}
\caption{Examples of gradient vector fields on cell complexes.
Each point in the interior of the edge $PQ$ has three directions.}\label{fig:1}
\end{figure}

\begin{defn}\label{def:gradfield}
The set of \emph{gradient vectors} of $f$ at $x$ is 
 the set of tangent vectors $\{\tau_{x;j, l}\}$.
The \emph{gradient vector field} of $f$ is the assignment of the set of gradient vectors
 to each point of $X_1$.
\end{defn}

This definition of gradient vector fields is very different from the usual gradient vector field in the 
 smooth case, see Figure \ref{fig:1}.
There will be different definitions of vector fields in piecewise linear theory, see \cite{LP, St} for example.
In these references, vector fields on manifolds with piecewise linear structures are studied.
The definition here is different from these definitions as well, in that it allows multivaluedness.
This makes it possible to develop a meaningful theory of gradient flows on spaces 
 which are not necessarily manifolds.
On the other hand, the definition here is very simple compared to the other definitions in that 
 it is constant valued on cells.
Of course, such a trivial choice will not give a meaningful theory in general.
It makes sense only if it is combined with 
 the tameness of the function and the sharpness of the metric.

A nice feature of this definition of gradient vector fields is that, unlike the case of 
 discrete gradient paths (Definition \ref{def:gradpath}),
 it is purely determined by 
 the function and the metric, without referring to other information 
 such as the dimension of the vertices.
From this, one can immediately see the following duality.
Assume that the affine complex $X$ is a polyhedral decomposition of a manifold without boundary.
Consider the dual of our piecewise linear theory on the dual complex $\bar X$ of $X$.
Namely, since there is a natural inclusion-reversing 
 one-to-one correspondence between the set of cells of $X$ and $\bar X$, 
 we have a natural function $\bar F$ on the set of 
 cells of $\bar X$ induced from the discrete Morse function $F$ on $X$.
Then, $-\bar F$ is a discrete Morse function on $\bar X$.
Note that if $F$ is tame and generic, so is $-\bar F$.
Also, there is a natural identification between $X_1$ and the barycentric subdivision of $\bar X$.
Thus, we have a piecewise linear theory on $X_1$ using the function induced from $-\bar F$ and
 the given sharp piecewise affine metric.
\begin{prop}\label{prop:dualtheory}
The associated gradient vector for the dual theory at a point $x\in X_1$
 is given by $\{-\tau_{x;j, l}\}$.
\end{prop}
\proof
This follows from the definition. \qed\\

This point is one of the advantages of considering piecewise linear theory,
 compared to the discrete theory.
That is, in the discrete case, when one considers
 Poincar\'{e} duality for a polyhedral decomposition of a manifold without boundary, 
 the dual theory is defined on the dual complex of the original, so that one cannot  
 directly relate the original discrete gradient path to the dual.
Here, we can consider Poincar\'{e} duality on the same complex just by
 reversing the direction of the gradient vector field.

\section{Properties of piecewise linear flows}\label{sec:6}

%

%
Proposition \ref{prop:lowestvalue} below will be a basic ingredient in the argument in the following sections.
First, we prove the following.

\begin{lem}\label{lem:triangleflow}
Consider a triangle $P$ with vertices $a, b, c$ and put an affine metric 
 on it which is not necessarily sharp.
Take an affine linear function $f$ on $P$ which is generic so that 
 the values at the vertices are different, say, $f(a)<f(b)<f(c)$, and the 
 gradient flow with respect to the metric is not parallel to any of the edges. 
Then, if the edge $\overline{bc}$ of $P$ 
 has an in-flow, the angle $\angle abc$ is larger than $\frac{\pi}{2}$.
Further, the edge $\overline{ab}$ also has an in-flow.
\end{lem}
\proof
There is a point $b'$ on the edge $\overline{ac}$
 with the property $f(b) = f(b')$.
Since $f$ is affine linear, $f$ is constant on the line connecting $b$ and $b'$.
This line divides $P$ into two triangles 
 $P_1 = \triangle{cbb'}$ and $P_2 = \triangle{bab'}$.
The gradient flow is perpendicular to $\overline{bb'}$, and as an edge of 
 $P_1$, clearly $\overline{bb'}$ has an in-flow.
By the assumption, the edge $\overline{bc}$ also has an in-flow.
It follows that there is a flow into the vertex $b$ 
 whose trajectory is perpendicular to $\overline{bb'}$.
Then, we have
\[
\angle abc = \angle abb'+ \angle cbb' > \angle abb' + \frac{\pi}{2} >\frac{\pi}{2}.
\] 
The final assertion is obvious.\qed

\begin{prop}\label{prop:lowestvalue}
Let $X, X_1, F$ and $f$ be as above.
Put a sharp piecewise affine metric on $X_1$.
Let $\tau$ be a cell in $X_1$. 
If $\tau$ has an in-flow from a cell $\sigma$, the cell $\tau$
 contains the lowest value vertex of the Morse function on $\sigma$.
\end{prop}
\proof
Note that $\tau$ is a facet of $\sigma$ by definition, see Definition \ref{def:in-out}.
If $\dim\sigma = 1$, the claim is obvious.
If $\dim\sigma = 2$, the claim follows from Lemma \ref{lem:triangleflow}.
Namely, using the notation in Lemma \ref{lem:triangleflow}, if the cell $\tau$
 does not contain the lowest value vertex, we have $\tau = \overline{bc}$.
Thus, if $\tau$ has an in-flow, the angle $\angle abc$ is larger then $\frac{\pi}{2}$.
However, this contradicts to the sharpness of the metric.

Assume $\dim\sigma$ is larger than two.
Suppose the lowest value vertex $p$ of $\sigma$ is not contained in $\tau$.
Take a generic point $q$ on $\tau$.
Let $L$ be a two dimensional plane containing the edge $\overline{pq}$ and
 the flow line from the interior of $\sigma$ to $q$.
Then, the intersection $\tau\cap L$ is a segment, which we write as $\overline{rs}$.
We can assume the inequalities
\[
f(p) < f(r) < f(s)
\]
 hold.
Then, the triangle $\triangle prs$ satisfies the condition of Lemma \ref{lem:triangleflow}, 
 since there is an in-flow on $\overline{rs}$ by the assumption.
It follows that the angle $\angle prs$ is larger than $\frac{\pi}{2}$.
However, this contradicts to the sharpness of the metric.\qed\\

Dually, we have the following.
\begin{cor}\label{cor:largestvalue}
Using the same notation as in Proposition \ref{prop:lowestvalue}, 
 if $\tau$ has an out-flow to a cell $\sigma$, the cell $\tau$
 contains the largest value vertex of the Morse function on $\sigma$.\qed
\end{cor}

\section{Gradient trajectories between critical vertices}\label{sec:7}
Let $X, X_1, F$ and $f$ be as before.
Fix a sharp piecewise affine metric on $X_1$.

%

%
%

%
\begin{defn}
Critical points of an affine complex $X$ equipped with 
 a discrete Morse function are the vertices of $X_1$ corresponding
 to the critical cells of the discrete Morse function in the
 sense of Forman. 
We call a critical point $i$-dimensional
 if it is a vertex corresponding to
 an $i$-dimensional cell of the original complex.
\end{defn}

As in discrete theory, let $\mathcal M_i^{PL}$ be the free abelian group generated by 
 critical $i$-dimensional vertices.
It is canonically isomorphic to $\mathcal M_i$.
We want to construct a complex 
\[
0\to \mathcal M_n^{PL}\to \cdots \to \mathcal M_0^{PL}\to 0
\]
 by counting gradient trajectories of the piecewise linear flow.
First, let us define gradient trajectories.

\begin{defn}
Let $\sigma$ be a cell of $X_1$.
Let $x, y$ be points on $\sigma$.
A \emph{gradient segment} from $x$ to $y$  
 is an affine linear map 
\[
\gamma\colon [t_0, t_1]\to \sigma,
\]
 where $t_0, t_1\in \Bbb R$, 
 such that $\gamma(t_0) = x, \gamma(t_1) = y$ and
 for each $s\in [t_0, t_1]$, $\gamma'(s)\in \{\tau_{x;j, l}\}$ (see Definition \ref{def:gradfield}).
A \emph{gradient trajectory} from $x$ to $y$ is a piecewise affine linear map 
\[
\gamma\colon [t_0, t_k]\to \sigma,
\]
 where $k$ is a positive integer, such that there is a refinement
\[
t_0 < t_1 < \cdots < t_k
\]
 of $[t_0, t_k]$ where $\gamma|_{[t_{l-1}, t_l]}$ is a gradient segment for
 each $l = 1, \dots, k$.
\end{defn}
\begin{defn}\label{def:gradtraj}
A gradient trajectory of the gradient vector field on $X_1$
 is an ordered sequence
\[
\gamma_1, \gamma_2, \dots, \gamma_k,
\]
 of gradient trajectories $\gamma_l: [t_{l-1}, t_l]\to \sigma_l$
 on some cell $\sigma_l$ such that 
\[
\gamma_l(t_l) = \gamma_{l+1}(t_l), \;\; l = 1, \dots, k-1.
\] 
\end{defn}
We identify two gradient trajectories if their images give the same
 subsets of $X$.
Note that for a point $x\in X_1$, there can be more than one gradient trajectories
 through $x$, contrary to the case of smooth manifolds. 
 
\section{Piecewise linear Morse complex}\label{sec:8}
\begin{defn}\label{def:rib}
Let $i$ be a non-negative integer.
Let $X^{(i)}$ be the $i$-th skeleton of $X$.
Let $X_1^{(i)}$ be the barycentric subdivision of $X^{(i)}$.
We call $X_1^{(i)}$ the \emph{$i$-th rib} of $X_1$.
\end{defn} 
Note that $X_1^{(i)}$ is not the $i$-th skeleton of $X_1$, 
 which consists of all the simplices of $X_1$ of dimension at most $i$.

\begin{prop}\label{prop:flow}
Assume $F$ is generic and tame.
Let $p$ and $q$ be $i$- and $(i-1)$-dimensional critical vertices, respectively.
Then, a gradient trajectory $\gamma: [0, t]\to X_1$
 connecting $p$ and $q$ is contained in a subcomplex
 composed of edges whose ends are $i$- and $(i-1)$-dimensional vertices.
\end{prop}
\proof
First, we observe the following.
\begin{claim}\label{claim:flow}
Let $x\in X_1^{(i)}$ be a point.
Let $\delta$ be a gradient trajectory starting from $x$.
Then, the inclusion $\delta\subset X_1^{(i+1)}$ holds.
\end{claim}
\proof
This is an immediate consequence of Corollary \ref{cor:largestvalue}
 and the tameness of $F$.\qed\\

Since $p$ is an $i$-dimensional critical vertex, 
 any of its adjacent $(i+1)$-dimensional vertices has larger value of $f$ than $f(p)$.
Then, again by Corollary \ref{cor:largestvalue}, 
 for sufficiently small $\varepsilon$, 
 $\gamma([0, \varepsilon))$ is contained in $X^{(i)}$.
Let $\sigma$ be the cell of $X_1$ with the property that
 $\gamma((0, \varepsilon))\subset int(\sigma)$ for sufficiently small $\varepsilon$,
 where $int(\sigma)$ is the interior of $\sigma$.
Note that $\sigma$ is uniquely determined by this property.
Let $0<t_1$ be the minimal real number with the property
 $\gamma(t_1)\in \partial\sigma$.
Let $\tau\subset \partial\sigma$ be the unique cell with the property that
 $\tau$ contains $\gamma(t_1)$ in its interior.
Then, since $p$ is a critical vertex, the inclusion $\tau\subset X_1^{(i-1)}$ holds.

Thus, by Claim \ref{claim:flow}, the image of the restriction $\gamma_{[t_1, t]}$ is contained in $X_1^{(i)}$.
It follows that we have $\gamma([0, t])\subset X_1^{(i)}$.
Since $q$ is a critical vertex, a positive dimensional
 cell in $X_1^{(i)}$ which has $q$ as
 the lowest value vertex with respect to $f$ must be an edge connecting $q$ and
 an adjacent $i$-dimensional vertex.
Thus, by Proposition \ref{prop:lowestvalue},
 the intersection of a small neighborhood of $q$ with $\gamma$
 must be a subset of an edge connecting $q$ and an $i$-dimensional vertex $q_1$.
Similarly, by the tameness of $F$, a positive dimensional 
 cell in $X_1^{(i)}$ which has $q_1$ as
 the lowest value vertex with respect to $f$ must be an edge connecting $q_1$
 and an $(i-1)$-dimensional vertex $q_2$ (which is unique, if any).
In this case, we have
\[
f(q_1) < f(q_2).
\]
Note that any $(i-2)$-dimensional vertex adjacent to $q_2$ is also adjacent to $q_1$.
If $r$ is an $(i-2)$-dimensional vertex adjacent to $q_1$, 
 we have $f(r) < f(q_1)$ by the tameness of $F$.
In particular, for any $(i-2)$-dimensional vertex $r$ adjacent to $q_2$, we have $f(r)<f(q_2)$.

Thus, again, a positive dimensional 
 cell in $X_1^{(i)}$ which has $q_2$ as
 the lowest value vertex with respect to $f$ must be an edge
 connecting $i$- and $(i-1)$-dimensional vertices.
Repeating this, we see that $\gamma$ is contained in the union of edges
 connecting $i$- and $(i-1)$-dimensional vertices.\qed\\

From this, the following is obvious.

\begin{cor}
In the notation of Proposition \ref{prop:flow},
 the number of gradient trajectories connecting $p$ and $q$ is finite.\qed
\end{cor}

\subsection{Sign of gradient trajectory and Morse differential}\label{subsec:sign2}

We define the linear map 
\[
d_{PL}: \mathcal M_i^{PL}\to \mathcal M_{i-1}^{PL}
\]
 by linearly extending the map
\[
p\mapsto \sum_{q\in \mathcal M_{i-1}^{PL}}\sum_{\gamma\in\Gamma_{PL}(p, q)}\bar m(\gamma) q, 
\]
 where $\Gamma_{PL}(p, q)$ is the set of gradient trajectories connecting $p$ and $q$, 
 and $\bar m(\gamma)$ is the sign of the gradient trajectory $\gamma$ determined as follows.
Namely, to each $k$-dimensional critical vertex $a^{(k)}$, $0\leq k\leq n$, we attach an orientation of the 
 $k$-dimensional cell
 corresponding to it.
It is represented by a constant element of $\Gamma(\sigma_{a^{(k)}}, \wedge^kT\sigma_{a^{(k)}})$.
As an element of $\mathcal M_k^{PL}$, $-a^{(k)}$ represents the same vertex with the opposite orientation
 of $\sigma_{a^{(k)}}$ attached.

By Proposition \ref{prop:flow}, 
 the trajectory connecting $p^{(i)}\in \mathcal M_i^{PL}$ and $q^{(i-1)}\in \mathcal M_{i-1}^{PL}$
 is contained in the union of edges connecting $i$- and $(i-1)$-dimensional vertices.
It follows that if the trajectory passes through the vertices as
\[
p^{(i)} \to s^{(i-1)}\to r^{(i)} \to \cdots \to q^{(i-1)}, 
\]
 the orientation at $p^{(i)}$ naturally induces an orientation at each $s^{(i-1)}, r^{(i)}, \dots, q^{(i-1)}$.
Namely, if $o_{p^{(i)}}\in \Gamma(\sigma_{p^{(i)}}, \wedge^iT\sigma_{p^{(i)}})$ is the orientation at $p^{(i)}$, 
 it induces an orientation at $s^{(i-1)}$ by
\[
o_{s^{(i-1)}} = n_{p^{(i)}}\lfloor o_{p^{(i)}} \in \Gamma(\sigma_{s^{(i-1)}}, \wedge^{(i-1)}T\sigma_{s^{(i-1)}}), 
\]
 where $n_{p^{(i)}}$ is the unit outward normal vector of
 $\sigma_{p^{(i)}}$ on the facet $\sigma_{s^{(i-1)}}$, and 
 $\lfloor$ denotes the interior product defined using the affine metric.
Then, the orientation at $r^{(i)}$ is given by
\[
o_{r^{(i)}} = -n_{r^{(i)}}\wedge o_{s^{(i-1)}}, 
\]
 where $n_{r^{(i)}}$ is the unit outward normal vector of
  $\sigma_{r^{(i)}}$ on the facet $\sigma_{s^{(i-1)}}$.
Repeating this, an orientation is induced at $q^{(i-1)}$.
If we compare this with the given orientation, we set $\bar m(\gamma) = 1$
 when they coincide, and $-1$ if not.

With this definition, we have a sequence of maps
\begin{equation}\label{eq:plcpx}
0\to \mathcal M_n^{PL}\xrightarrow{d_{PL}}\mathcal M_{n-1}^{PL}
 \xrightarrow{d_{PL}}\cdots \xrightarrow{d_{PL}} \mathcal M_0^{PL}\to 0.
\end{equation}
A priori, we do not know whether this gives a chain complex.
In the next section, we see this is indeed the case, and its homology
 is the same as the singular homology of $X$.
\section{The main theorem}\label{sec:9}

Let $X, X_1, F$ and $f$ be as above.
Namely, $X$ is an affine complex of dimension $n$, 
 $X_1$ is the barycentric subdivision of $X$, 
 $F$ is a tame generic discrete Morse function on $X$, and
 $f$ is the affine linear function on $X_1$ constructed from $F$.
Put a sharp piecewise affine metric on $X_1$.
Our purpose is to show the following result.
\begin{thm}\label{thm:main}
The sequence $({\mathcal M}_{\bullet}^{PL}, d_{PL})$ in the previous section
 is a complex, and its homology is isomorphic to
 the singular homology of $X$.
\end{thm} 
\proof

We examine the relation between discrete and piecewise linear Morse theories.
Recall that a gradient path in discrete Morse theory (Definition \ref{def:gradpath})
 is an appropriate sequence $\gamma$ of $i$-cells
 of $X$
\begin{equation*}
\gamma = \alpha_0, \alpha_1, \alpha_2, \dots, \alpha_r.
\end{equation*}
Let $\Gamma_r(\alpha_0, \alpha_r)$ be the set of such gradient paths.

Consecutive $i$-cells $\alpha_{l-1}$ and $\alpha_l$ are in the boundary of a uniquely determined
 $(i+1)$-cell, if $\alpha_{l-1}\neq \alpha_l$.
Let us write this $(i+1)$-dimensional cell by $\beta_{l}$.
If $\alpha_{l-1} = \alpha_l$, take $\beta_{l} = \emptyset$.
Let  $\beta_0$ be an $(i+1)$-dimensional cell adjacent to $\alpha_0$ with the property 
 $F(\beta_0)>F(\alpha_0)$, 
 and consider the sequence
\[
\beta_0, \alpha_0, \beta_1, \alpha_1, \beta_2, \cdots, \alpha_{r-1}, \beta_r, \alpha_r.
\]
Let $\Gamma_r(\beta_0, \alpha_r)$ denote the set of such sequences.
An element of $\Gamma_r(\beta_0, \alpha_r)$
 is determined by an element of $\Gamma_r(\alpha_0, \alpha_r)$ and
 the choice of $\beta_0$.
Then, if $\beta_0$ is an $(i+1)$-dimensional critical cell,
  the differential of discrete Morse theory (\ref{eq:discretediff})
 can be written in the form
\begin{equation}
\tilde \del \beta_0 = \sum_{\alpha\in\mathcal M_i}
 \sum_{\tilde{\gamma} \in \Gamma_N(\beta_0, \alpha)} \tilde m(\tilde{\gamma})\alpha,
\end{equation}
 where $N$ is a sufficiently large number and
 $\mathcal M_i$ is the set of critical $i$-cells.
The sign
 $\tilde m(\gamma)\in\{\pm 1\}$ is equal to $m(\gamma)$, where $\gamma$ is the element of 
 $\Gamma_N(\alpha_0, \alpha)$ associated with $\tilde{\gamma}$.
Here, $\alpha_0$ is a facet of $\beta_0$ with an induced orientation,
 that is, $\langle \partial\beta_0, \alpha_0\rangle = 1$.


To adjust the terminology with those in piecewise linear theory, we redefine the sequence 
\[
\beta_0, \alpha_0, \beta_1, \alpha_1, \beta_2, \cdots, \alpha_{r-1}, \beta_r, \alpha_r
\]
 above 
 as a gradient path in discrete Morse theory.
We identify two sequences $\{\beta_i, \alpha_i\}$ and $\{\beta_i', \alpha_i'\}$
 which are the same up to stabilization, 
 that is, if the union $\cup_l int({\beta}_l)\bigcup\cup_l\alpha_l$
 and $\cup_l int({\beta}'_l)\bigcup\cup_l\alpha'_l$ are the same as a subset of $X$.
Here, $int({\beta}_l)$ and $int({\beta}'_l)$ are the interior of $\beta_l$ and $\beta_l'$, respectively.

Clearly, an element in $\Gamma_N(\beta_0, \alpha_r)$
 determines a sequence
\[
p_0\to q_0\to p_1\to\cdots \to p_r\to q_r,
\]
 where $p_i$ is the $(i+1)$-dimensional vertex in $X_1$
 corresponding to $\beta_i$,
 $q_i$ is the $i$-dimensional vertex in $X_1$ corresponding to $\alpha_i$, 
 and $r$ is the smallest number at which the sequence
 $\alpha_0, \alpha_1, \alpha_2, \dots, \alpha_r$
 stabilizes (that is, the smallest number such that $\alpha_r =\alpha_{r+1}$ holds).
This sequence has the following properties:
\begin{itemize}
\item 
Consecutive vertices in this sequence are connected by 
 a unique edge in $X_1$.
\item Following inequalities hold:
\[
f(p_0)>f(q_0)>f(p_{1}) > \cdots > f(p_r) > f(q_r).
\]
\end{itemize}
Conversely, a sequence of $i$- and $(i+1)$-dimensional adjacent vertices
 with the properties above 
 determines a unique (up to stabilization) gradient path in discrete Morse theory.
By Proposition \ref{prop:flow}, we have the following.
%
\begin{prop}\label{prop:corresp}
Let $\beta$ be a critical $(i+1)$-dimensional cell and $\alpha$ be a
 critical $i$-dimensional cell
 of discrete Morse theory.
Let $p, q$ be corresponding $(i+1)$- and $i$-dimensional critical vertices of 
 piecewise linear Morse theory.
There is a canonical one to one correspondence between the following objects:
\begin{itemize}
\item Gradient paths in 
 discrete Morse theory connecting $\beta$ and $\alpha$ up to stabilization. 
\item Gradient trajectories in piecewise linear Morse theory connecting $p$ and $q$. 
\end{itemize}
\end{prop}

Recall that the differential $d_{PL}\colon \mathcal M_{i+1}^{PL}\to \mathcal M_i^{PL}$
 in piecewise Morse theory is given by 
\[
d_{PL}p = \sum_{q\in \mathcal M_{i}^{PL}}\sum_{\gamma\in\Gamma_{PL}(p, q)}\bar m(\gamma) q,
\]
 and 
 the differential $\tilde{\partial}\colon \mathcal M_{i+1}\to \mathcal M_i$
   in discrete Morse theory is given by 
\[
\tilde \del \beta = \sum_{\alpha\in\mathcal M_i}
  \sum_{\tilde{\gamma} \in \Gamma_N(\beta, \alpha)} \tilde m(\tilde{\gamma})\alpha.
\]
Let $\tilde{\gamma}$ be the element of $\Gamma_N(\beta, \alpha)$
 corresponding to $\gamma\in \Gamma_{PL}(p, q)$ by Proposition \ref{prop:corresp}. 
Then, Theorem \ref{thm:main} is a consequence of the following claim and Theorem \ref{thm:discmtheory}.
 
\begin{prop}\label{prop:sign}
The equality $\bar m(\gamma) = \tilde m(\tilde{\gamma})$ holds.
\end{prop}
\proof
This follows from a straightforward comparison between the argument in 
 Subsections \ref{subsec:sign1} and \ref{subsec:sign2}.\qed\\

\section{Stable and unstable complexes}\label{sec:stablecomplex}
Let $X, X_1, F$ and $f$ be as before.
Let $i$ be a positive integer.
In this section, we will show that in piecewise linear Morse theory, we have nice stable and 
 unstable complexes associated with critical points.
This is another benefit from the fact that we can consider the original Morse theory and its dual 
 in the same complex, see Proposition \ref{prop:dualtheory}.
We use the same notation as in the previous sections.
First, we observe the following.
\begin{lem}
Let $C$ be a closed subset of $X_1$.
Then, the union of gradient trajectories starting from a point in $C$ is a closed subset of $X_1$.
\end{lem}
\proof
It suffices to prove the claim on a simplex $\sigma$.
If $x_1$ is a point in the interior of $\sigma$, there is a unique
 gradient vector at $x_1$ tangent to $\sigma$.
This does not depend on the choice of $x_1$, and gradient vectors at the points on the boundary
 of $\sigma$ also contain this direction.
It follows that the subset of $\sigma$ swept by the flow along this direction starting from 
 points in $C$ is closed in $\sigma$.
Namely, embed $\sigma$ in some $\Bbb R^d$ isometrically and let $v\in\Bbb R^d$
 be the direction of the gradient flow on $\sigma$.
If we define the set $C_{\geq}$ by
\[
C_{\geq} = \{x+tv\;|\; x\in C,\; t\geq 0\}\subset \Bbb R^d, 
\]
 the subset of $\sigma$ swept by the flow starting from points in $C$
 is $C_{\geq}\cap \sigma$.
Since $C_{\geq}$ is closed, the claim follows.

The same argument applies to the intersection of a face of $\sigma$ and $C$.
Since the union of gradient trajectories starting from a point in $C$ is
 the union of these closed subsets on each face of $\sigma$, the claim follows.\qed
 
\begin{lem}\label{lem:traj}
Let $\sigma$ be a simplex in $X_1$.
Let $p$ be the vertex of $\sigma$ at which the function $f$ takes the highest value in $\sigma$.
Then, $\sigma$ is swept by the gradient trajectories starting from $p$.
\end{lem}
\proof 
By the previous lemma, it suffices to prove that the interior of $\sigma$
 is swept by the gradient trajectories starting from $p$.
We use an induction.
If $\dim\sigma =1$, the claim is obvious.
Assume that the claim is proved for $\sigma$ with $\dim\sigma\leq k$
 for some positive integer $k$.
Let us prove the case where $\dim\sigma = k+1$.

Let $x$ be a point in the interior of $\sigma$.
As in the proof of the previous lemma, there is a unique gradient vector at $x$ 
 tangent to $\sigma$.
Tracing back the gradient trajectory along this direction, we will hit a boundary of 
 $\sigma$.
Let $q$ be the point in $\partial\sigma$ at which this occurs for the first time.
Let $\{A_1, \dots, A_m\}$ be the set of facets of $\sigma$ containing $q$.
By Corollary \ref{cor:largestvalue}, each $A_l$ contains $p$.
Then, by applying the induction hypothesis to any of the cells $\{A_1, \dots, A_m\}$, 
 the claim follows.\qed.

\begin{lem}\label{lem:sweep}
Let $q$ be an $(i-1)$-dimensional vertex and $p$ be an adjacent $i$-dimensional vertex
 such that $f(q)>f(p)$ holds.
Let $\sigma_p$ be the cell of $X$ corresponding to $p$.
Then, $\sigma_p$ is swept by the gradient trajectories starting from $q$. 
\end{lem}
\proof
By the tameness, if $\alpha$ is a maximal dimensional
 simplex in $\sigma_p$, the vertex taking the highest value among the 
 set of vertices of $\alpha$ is $p$ or $q$.
Since $p$ is contained in a gradient trajectory starting from $q$, the claim follows 
 from the previous lemma.\qed\\

We also prove the following claim.
\begin{lem}\label{lem:confine}
Let $p$ be an $i$-dimensional critical vertex.
Let $\gamma$ be a gradient trajectory starting from $p$.
Then, $\gamma$ is contained in $X_1^{(i)}$, the $i$-th rib of $X_1$ (see Definition \ref{def:rib}).
\end{lem}
\proof
We assume $i>0$ since otherwise the claim is trivial.
By Corollary \ref{cor:largestvalue}, for sufficiently small $\varepsilon$, the image 
 $\gamma((0, \varepsilon))$ is contained in $X_1^{(i)}$.
Let $\tau$ be the unique cell of $X_1$ satisfying $\gamma((0, \varepsilon))\subset int(\tau)$.
Let $t_1>0$ be the minimal real number satisfying $\gamma(t_1)\in \partial\tau$.
Then, we have $\gamma(t_1)\in X_1^{(i-1)}$.
The lemma follows from this and Claim \ref{claim:flow}.\qed\\

Let $k$ be a positive integer.
\begin{prop}\label{prop:subcomplex}
Let $p$ be an $i$-dimensional critical vertex.
Then, the subset of $X_1$ swept by the trajectories starting from $p$
 is a subcomplex of $X^{(i)}$.
\end{prop}
\proof
%
By Corollary \ref{cor:largestvalue}, for sufficiently small $\varepsilon$, the image 
 $\gamma((0, \varepsilon))$ is contained in $\sigma_p$, for any trajectory $\gamma$ 
 starting from $p$.
By Lemma \ref{lem:traj}, the cell $\sigma_p$ is swept by the trajectories starting from $p$.
Thus, to prove that the union of the trajectories is a subcomplex of $X$, it suffices to show that 
 for any facet $A$ of $\sigma_p$, the union of trajectories starting from points on $A$
 is a subcomplex of $X$.
If this is proved, by Claim \ref{claim:flow}, it must be a subcomplex of $X^{(i)}$.

We prove this from the following more general statement.
\begin{claim}
Let $B$ be a cell in $X$.
Then, the union of trajectories starting from points on $B$
 is a subcomplex of $X$.
\end{claim}
\proof
We prove this by induction on the dimension of $B$.
First, let us consider the case when $\dim B = 0$. 
If $B$ is a critical cell, the claim is obvious since there is no trajectory emanating from $B$. 
If $B$ is not critical, there is a unique one dimensional cell $C$ adjacent to $B$ such that 
\[
F(B)>F(C)
\]
 holds.
Higher dimensional cells adjacent to $B$ have larger value of $F$ than $F(B)$ by the tameness. 
 
By Corollary \ref{cor:largestvalue}, any trajectory starting from $B$ is given by an out-flow
 into $C$.
On the other hand, again by the tameness, there is only one cell $B'$ adjacent to $C$
 satisfying 
\[
F(C)>F(B'),
\]
 and this $B'$ is a zero dimensional cell.
Repeating this argument, the claim follows.
 
Assume that we have proved the claim when $\dim B\leq k-1$ for some positive integer $k$.
Let us prove the claim when $\dim B = k$.

We use further induction on the value of $F$.
Let $B_1$ be a $k$-dimensional cell of $X$ such that there is no other $k$-dimensional cell
 whose value of $F$ is smaller than $F(B_1)$.
Then, there is no adjacent $(k+1)$-cell which has a smaller value of $F$ than $F(B_1)$, 
 by definition of discrete Morse functions.
Then, $B_1$ is critical, or there is a facet $C$ of $B_1$ satisfying $F(C) > F(B_1)$.

If $B_1$ is critical, as in the beginning part of the proof of Proposition \ref{prop:subcomplex}, 
 any trajectory starting from a point on $B_1$ touches $\partial B_1$.
Thus, the claim follows by the induction hypothesis on the dimension of $B$.
If there is a facet $C$ of $B_1$ satisfying $F(C) > F(B_1)$, $B_1$ is swept by 
 trajectories starting from $C$, by Lemma \ref{lem:sweep}.
Thus, it suffices to show that the union of trajectories starting from $C$ is a subcomplex of $X$, 
 since it coincides with the union of trajectories starting from $B_1$.
However, this follows from the induction hypothesis on dimension again.
 
Now, let 
\[
c_1<c_2<\cdots < c_l
\]
 be the values of $F$ taken by $k$-dimensional cells of $X$.
Assume that the claim is proved for a $k$-dimensional cell whose value of $F$ is at most $c_i$, 
 $1\leq i<l$.
We will prove the claim for a $k$-dimensional cell $B$ satisfying $F(B) = c_{i+1}$.

If $B$ is critical or there is a facet $C'$ of $B$ satisfying $F(C') > F(B)$, 
 the claim is proved as in the case of $B_1$ above.
Assume that there is a $(k+1)$-dimensional cell $D$ adjacent to $B$ satisfying $F(B)>F(D)$.
By Lemma \ref{lem:sweep}, $D$ is swept by trajectories starting from $B$.
Moreover, if $x\in int(B)$, for sufficiently small $\varepsilon$, 
 a trajectory $\gamma$ starting from $x$ satisfies $\gamma((0, \varepsilon))\subset D$, 
 by definition of discrete Morse functions and Corollary \ref{cor:largestvalue}.
It follows that the trajectory $\gamma$ hits some facet of $D$ other than $B$ itself.
Therefore, it suffices to prove that the union of trajectories starting from a facet $E$
 of $D$ other than $B$ is a subcomplex of $X$.
However, we have $F(E) < F(B)$ by definition of discrete Morse functions.
Therefore, we can apply the induction hypothesis on the values of $F$ to $E$.
This proves the claim and the proposition.\qed\\

Using this proposition, we can define stable and unstable complexes of critical points in
 piecewise linear Morse theory.

\begin{defn}
Assume $X$ is a manifold without boundary.
We call the subcomplex of $X^{(i)}$ constructed in Proposition \ref{prop:subcomplex}
 the \emph{unstable complex} of the critical point $p$.
Consider the dual theory, see Section \ref{sec:5}.
Then, $p$ is an $(n-i)$-dimensional critical point.
We call its unstable complex the \emph{stable complex} of $p$
 of the original theory.
\end{defn}

Note that the stable complex is not a subcomplex of $X$, but of its dual complex.
In any case, it is a subcomplex of $X_1$.

Since the gradient trajectory is defined purely in terms of a function and a metric as we mentioned in
 Section \ref{sec:5},
 this definition can be extended to the case when $X$ has boundary, or
 is not necessarily a manifold.
Namely, if $p$ is a critical point, we define the unstable complex as the subcomplex of $X^{(i)}$ 
 as in Proposition \ref{prop:subcomplex}, and 
 define the stable complex as the subcomplex with respect to the reversed gradient flow.
Note that the reversed gradient flow is not necessarily induced by a discrete Morse function, 
 since if $X$ is an affine complex which is not a manifold, a dual complex does not make sense in general.
In other words, it is not a gradient flow associated with a discrete Morse function, although it is the 
 gradient flow associated with a piecewise linear function on $X_1$. 
Nevertheless, the set of gradient trajectories gives a complex, which is isomorphic to the dual of the original complex
 $\mathcal M$.



\section*{Acknowledgement}
\noindent
The author was supported by JSPS KAKENHI Grant Number 18K03313.

\end{document}